\documentclass[twocolumn,pra,aps,superscriptaddress]{revtex4}

\usepackage{color}
\usepackage{epsfig}
\usepackage{latexsym}
\usepackage{amssymb}
\usepackage{amsmath}
\usepackage{algorithm}
\usepackage{algorithmic}
\usepackage{wrapfig}
\usepackage[apple mac]{inputenc}
\usepackage[english]{babel}
\usepackage{times}
\usepackage{latexsym}
\usepackage{fancyhdr}
\usepackage{verbatim}
\usepackage{tabularx}
\usepackage{epsfig}
\usepackage{amsmath}
\usepackage{amssymb}
\usepackage{graphicx}
\usepackage{wasysym}

\usepackage{yfonts}

\newcommand{\qed}{\hspace*{\fill}$\square$}

\newcommand{\be}{\begin{equation}}
\newcommand{\ee}{\end{equation}}
 \newcommand{\N}{\mathbf{N}}

\begin{document}

\title{Chiral Prime Concatenations}
\author{Miguel A. Martin-Delgado}
\affiliation{Departamento de F\'{\i}sica Te\'orica, Universidad Complutense, 28040 Madrid, Spain.\\
CCS-Center for Computational Simulation, Campus de Montegancedo UPM, 28660 Boadilla del Monte, Madrid, Spain.}

\begin{abstract} 
The notion of chiral prime concatenations is studied as a recursive construction of prime numbers 
starting from a seed set and with appropriate blocks to define the primality growth, generation by generation, either from the right or 
from the left. Several basic questions are addressed like whether chiral concatenation is a symmetric process, an endless process, as well as the calculation of largest chiral prime numbers. In particular, the largest left-concatenated prime number is constructed. It is a unique prime number of 24 digits. By introducing anomalous left-concatenations of primes we can surpass the limit of 24 digits for left-concatenated primes. It is conjectured that prime numbers are left chiral under anomalous concatenations.
\end{abstract}

\maketitle

\tableofcontents

%%%%%%%%%%%%%%%%%%%%%%%%%%%%%%%%%%%%%%
%%%%%%%%%%%%%%%%%%%%%%%%%%%%%%%%%%%%%%
\section{Introduction}
\label{sec:introduction}
%%%%%%%%%%%%%%%%%%%%%%%%%%%%%%%%%%%%%%
%%%%%%%%%%%%%%%%%%%%%%%%%%%%%%%%%%%%%%

Integer sequences of prime numbers abound \cite{Weisstein2}. One of the  most important questions for them is whether they are infinite sequences of prime numbers. The infinitude of the primes was proved by Euclid \cite{Euclid} in a beautiful theorem. However,
the infinitude of the sequence of twin prime numbers is one of the major open problems in number theory \cite{Guy}. We study here two series of prime numbers constructed with the method of concatenation. Interestingly, the issue of their finiteness can be solved explicitly both for left and right concatenations.

A word of caution is in order: concatenation of prime numbers has nothing to do with super-prime numbers, also called higher-order primes that corresponds to prime-indexed primes (PIPs) \cite{superprime,superprime2}.  PIPs are primes forming a subsequence of the prime series with primes placed at prime positions. That is, if $p(i)$ is the ith prime, then a super-prime is of the form $p(p(i))$. This procedure can be iterated to produce even higher-order of prime numbers \cite{superprime3}. Terminology is important since prime concatenations were introduced originally under the name of superprimes \cite{quantum}. We further elaborate on the implications of prime concatenations.

%%%%%%%%%%%%%%%%%%%%%%%%%%%%%%%%%%%%%%
%%%%%%%%%%%%%%%%%%%%%%%%%%%%%%%%%%%%%%
\section{Prime Number Primer}
\label{sec:primer}
%%%%%%%%%%%%%%%%%%%%%%%%%%%%%%%%%%%%%%
%%%%%%%%%%%%%%%%%%%%%%%%%%%%%%%%%%%%%%

An integer $n\in \N$ can be constructed via addition by starting from the unit element 1 and summing 
as many 1s up to the desired number $n$,
\begin{equation}\label{integer}
n:=\sum_{i=1}^n 1
\end{equation}
Trivial as this may sound, there is an important point in this construction that will be of interest later on.
Namely, as addition is a commutative operation, the construction of the integer number $n$ can proceed 
by adding 1s to the initial 1, starting either from the left or from the right. The resulting integer $n$ is indistinguishable
from which procedure we use, left or right.

The notion of prime number comes into play once the integers are constructed. Prime numbers are the elementary integer 
numbers with respect to multiplication: $n$ is prime iff its only divisors are 1 and itself. Then, the fundamental theorem of arithmetic
allows us to decompose any integer into prime factors,
\begin{equation}
  n= \prod_{i=1,\ldots, k }   p_i^{n_i} =
  p_1^{n_1} p_2^{n_2} \cdots p_k^{n_k},
\end{equation}
where $p_1<p_2<\ldots <p_k$  are prime numbers and the exponents $n_i$   are positive integers. It is assumed by convention that the empty product is equal to 1 (the empty product corresponds to $k = 0$).

 In this sense, the growth of prime numbers under addition or multiplication is symmetrically distributed to the left or the right. In other words,
 there is no chirality (preferred direction).
 
 %%%%%%%%%%%%%%%%%%%%%%%%%%%%%%%%%%%%%%
%%%%%%%%%%%%%%%%%%%%%%%%%%%%%%%%%%%%%%
\section{Concatenated Prime Numbers}
\label{sec:concatenated}
%%%%%%%%%%%%%%%%%%%%%%%%%%%%%%%%%%%%%%
%%%%%%%%%%%%%%%%%%%%%%%%%%%%%%%%%%%%%%

 Let us introduce another way to construct prime numbers called concatenation. Given two integers $n$ and $m$, the concatenation
 of them, denoted as $C[n,m]$, is another integer given by the numerical expression of $\text{nm}$. 
 \begin{equation}
C[n,m]:=\text{nm}.
\end{equation}
 This is not to be confused by the multiplication $n\cdot m$. For example, if we take the first two primes 2 and 3, the resulting concatenation 
 \begin{equation}
C[2,3]:=23,
\end{equation}
is another prime integer, whereas if we concatenate 2 and 5
\begin{equation}
C[2,5]:=25,
\end{equation}
is no longer a prime.

Now, we can think of concatenation as a tool to construct new primes much as we do with the addition  \eqref{integer}. Instead of using the unit 1 as a seed, we use the four single-digit prime numbers 
\begin{equation}\label{seed}
{\cal S} := \{ 2,3,5,7 \}.
\end{equation}
Concatenation is not a commutative operation unlike addition. Thus, we can construct primes by left-concatenation or by right-concatenation. Let us call them Chiral Prime Concatenations. This notion is base dependent. We shall be concerned with primes in decimal base. The resulting list of primes will be different. We consider first right-concatenation.

%%%%%%%%%%%%%%%%%%%%%%%%%%%%%%%%%%%%%%
%%%%%%%%%%%%%%%%%%%%%%%%%%%%%%%%%%%%%%
\subsection{Right-Concatenated Prime Numbers}
\label{sec:concatenatedR}
%%%%%%%%%%%%%%%%%%%%%%%%%%%%%%%%%%%%%%
%%%%%%%%%%%%%%%%%%%%%%%%%%%%%%%%%%%%%%

\begin{table}
\begin{tabular}{ |l|l| }
  \hline
  \multicolumn{2}{|c|}{{\bf Algorithm} Prime Right-Concatenation} \\
  \hline
  1 & {\bf Set} N: generation number \\
  2 & {\bf Set}t Initialization S=\{2,3,5,7\} \\
  3 & {\bf Set} Blocks BR=\{1,3,7,9\} \\
  4 & {\bf Set} CR[1] = S \\
  5 & {\bf While} $k\leq N$ \ {\bf do}\\
  6 & {\bf Concatenate} BR to the right of CR[k] \\
  7 & {\bf Apply} Primality Test to CR[k] \\
  8 & {\bf Update} CR[k] list of right-concatenated prime numbers \\
  9 & CR[k+1] $\leftarrow$ CR[k]  \ {\bf end do}\\
  10 & {\bf Output} CR[N] \\
  \hline
\end{tabular}
\caption{Pseudocode for generating right-concatenated prime numbers from the seed \eqref{seed} and the building blocks \eqref{right_block}.}
\label{tableR}
\end{table}

From the first generation of prime numbers given by \eqref{seed}, we can construct the second generation by concatenating
 to the right of ${\cal S}$ all possible one-digit odd integers excluding the 5. Denote this set of right-concatenating blocks as 
\begin{equation}\label{right_block}
{\cal B}_R := \{ 1,3,7,9 \}.
\end{equation}
Excluding the 5 avoids composite integers in the concatenation, but nevertheless it is not guaranteed that the right-concatenated integer constructed in this way is prime. Thus, upon right-concatenation of ${\cal S}$ with ${\cal B}_R$, we obtain the following list,
\begin{equation}\label{}
\begin{matrix}
\{\{21, 23, 27, 29\}, \{31, 33, 37, 39\}, \\ 
\{51, 53, 57, 59\}, \{71, 73, 77, 79\}\}.
\end{matrix}
\end{equation}
Selecting the prime numbers from this list we arrive at the second generation of right-concatenated prime numbers, denoted as
\begin{equation}\label{}
{\cal C}_R[2]:= \{23, 29, 31, 37, 53, 59, 71, 73, 79\}.
\end{equation}
The third generation of right-concatenated primes is obtained by right-concatenation of ${\cal C}_R[2]$ with ${\cal B}_R$ and discarding the resulting non-prime integers. The result is
\begin{equation}\label{}
\begin{aligned}
{\cal C}_R[3]:= &
\{233, 239, 293, 311, 313, 317, 373, 379, \\ 
&593, 599, 719, 733, 739, 797\}.
\end{aligned}
\end{equation}
And so on and so forth. A natural question arises: is the list of right-concatenated prime numbers endless, i.e., infinite?
It is easy to program a first order recursion relation based on the seed  ${\cal S}$ and the blocks ${\cal B}_R$ in order to produce all the generations sequentially. See Table \ref{tableR} for a pseudocode of the algorithm. Interestingly enough, the sequence ends after the $8^{\text th}$ generation. Thus, there are no right-concatenated prime numbers of nine digits or bigger.

Notice that all digits in the right-concatenated primes are odd numbers except for the digit 2 appearing as a seed (see App.\ref{sec:right})
to avoid even divisors. This is a consistency check of the algorithm. This is in sharp contrast with left-concatenated prime numbers
for which digits with even numbers abound (see Sec.\ref{sec:concatenatedL}).

\begin{figure}[t]
  \includegraphics[width=0.5\textwidth]{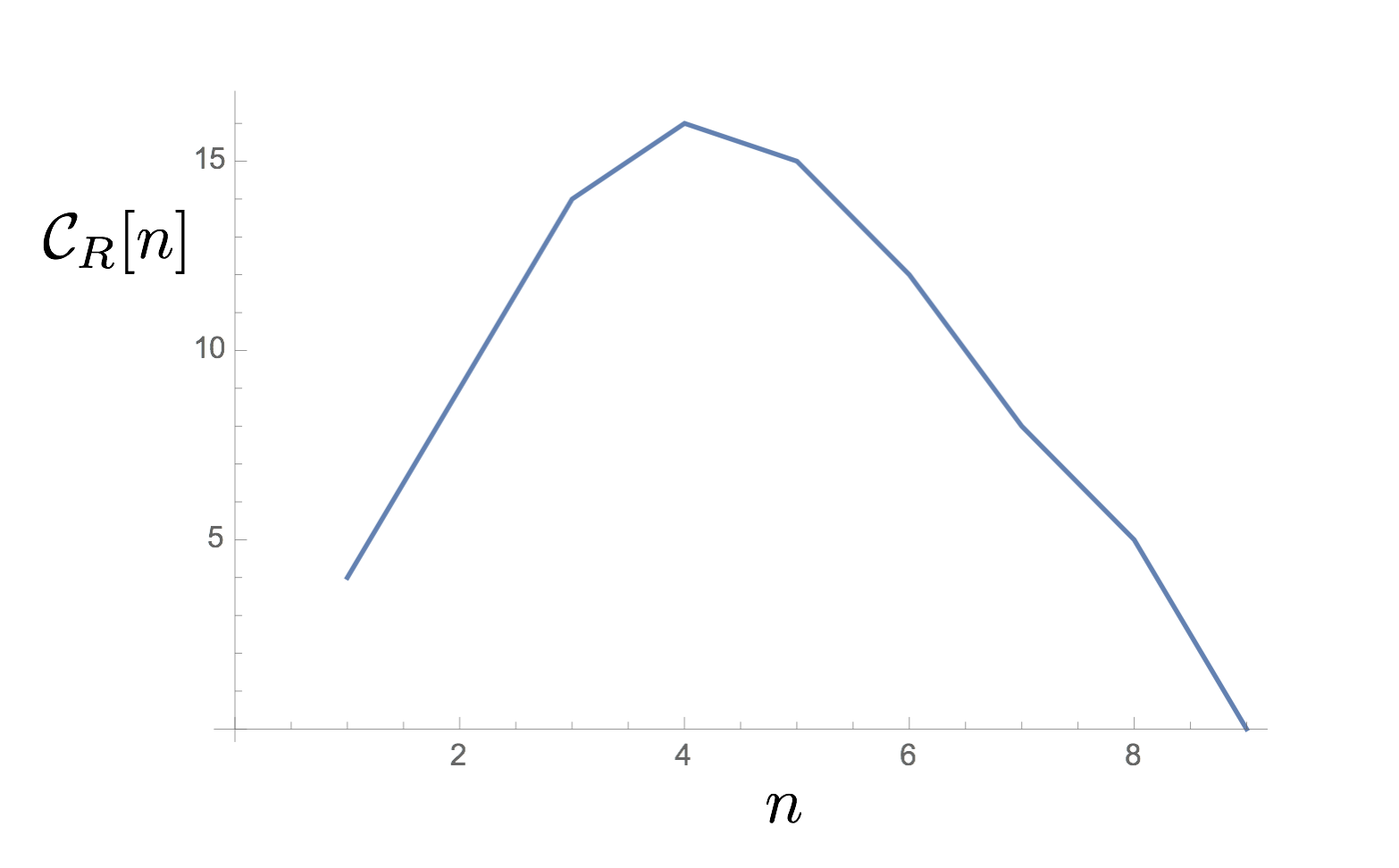}
 \caption{ The number of primes obtained by right-concatenation ${\cal C}_R[n]$ using the algorithm of Table \ref{tableR} as 
 a function of the generation number $n$. After the generation $n=8$ the right-concatenation procedure stops \cite{quantum}.}
  \label{fig:right}
\end{figure}

It is interesting to plot the number of right-concatenated prime numbers per generation. This is shown in Fig.\ref{fig:right}. After reaching a maximum, it falls down to zero. Notice that the largest right-concatenated prime number is not unique (see App.\ref{sec:right}). The behaviour of this plot is discussed in Sec.\ref{sec:unsolved}.

%%%%%%%%%%%%%%%%%%%%%%%%%%%%%%%%%%%%%%
%%%%%%%%%%%%%%%%%%%%%%%%%%%%%%%%%%%%%%
\subsection{Left-Concatenated Prime Numbers}
\label{sec:concatenatedL}
%%%%%%%%%%%%%%%%%%%%%%%%%%%%%%%%%%%%%%
%%%%%%%%%%%%%%%%%%%%%%%%%%%%%%%%%%%%%%

The construction of left-concatenated prime numbers starts with the same set of seeds \eqref{seed}, but now the
set of blocks is enlarged to include any single digit integer:
\begin{equation}\label{left_block}
{\cal B}_L := \{ 1,2,3,4,5,6,7,8,9 \}.
\end{equation}
Notice that we can remove the 2 from the seed since after the first generation, it will produce no primes because the only even prime number is 2. The reason to avoid the 0 in ${\cal B}_L $ is clearly different than in the right-concatenations where it gives divisible numbers. 
Firstly, we want to define concatenation by a single-step recursion relation in order to have a well-defined control over its growth and to be able to find concrete results unambiguously.
Secondly, we remove the 0 in order to have digits in ${\cal B}_L $ that are on equal footing among themselves. If we allow the 0 in the growth process from the left then we may add arbitrarily large numbers of 0s and then add some non-zero digit producing a prime number. To avoid these gaps of 0s we add the constraint of removing the 0 from ${\cal B}_L $. In Subsec.\ref{sec:anomalous}. 
we consider this case explicitly.

\begin{table}
\begin{tabular}{ |l|l| }
  \hline
  \multicolumn{2}{|c|}{{\bf Algorithm} Prime Left-Concatenation} \\
  \hline
  1 & {\bf Set} N: generation number \\
  2 & {\bf Set}t Initialization S=\{2,3,5,7\} \\
  3 & {\bf Set} Blocks BL=\{1,2,3,4,5,6,7,8,9\} \\
  4 & {\bf Set} CL[1] = S \\
  5 & {\bf While} $k\leq N$ \ {\bf do}\\
  6 & {\bf Concatenate} BL to the left of CL[k] \\
  7 & {\bf Apply} Primality Test to CL[k] \\
  8 & {\bf Update} CL[k] list of right-concatenated prime numbers \\
  9 & CL[k+1] $\leftarrow$ CL[k] \ {\bf end do}\\
  10 & {\bf Output} CL[N] \\
  \hline
\end{tabular}
\caption{Pseudocode for generating left-concatenated prime numbers from the seed \eqref{seed} and the building blocks \eqref{left_block}.}
\label{tableL}
\end{table}

The second generation of left-concatenated prime numbers is obtained by concatenating ${\cal B}_L $ to the left of the seed
${\cal S}$. Thus,
\begin{equation}\label{left-second}
{\cal C}_L[2] := \{ 13, 17, 23, 37, 43, 47, 53, 67, 73, 83, 97\}.
\end{equation}
Similarly, we can program a first order recursion relation starting from the seed ${\cal B}_L $ and concatenate to the left with the blocks of ${\cal B}_L$. See Table \ref{tableL} for a pseudocode of the algorithm.

The largest left-concatenated prime is unique and appears at the $24^{\text th}$ generation:
\begin{equation}\label{largest}
{\cal C}_L[24]:=
\{357686312646216567629137\}.
\end{equation}
It is generated from the prime number 7. By removing the digits from left to right, we obtain the largest
sequence possible of prime numbers. The resulting integer prime sequence is shown in App.\ref{sec:sequence}.
In the prime number ${\cal C}_L[24]$, half of the digits are odd numbers
and half are even numbers. Among the even digits, the 6 is the one that appears the most: seven times. Next comes
the 2 with three times and the 4 and the 8 once apiece. Among the odd digits, the 1, the 3 and the 7 appear three times each,
whereas the 5 appears twice and the 9 only once. The largest subsequence of even numbers has length five: 26462. 
By removing digits from the right we obtain non-prime numbers until arriving to the following sequence of seven digits 
3576863 that turns out to be prime again.
Interestingly enough, left-concatenated prime numbers at generations
$23^{\text rd}$, $22^{\text nd}$ and $21^{\text st}$ are also  generated from the 7, though not unique, namely,
\begin{equation}\label{}
\begin{aligned}
{\cal C}_L[23]:=&
\{57686312646216567629137, \\
&95918918997653319693967, \\
&96686312646216567629137\}.
\end{aligned}
\end{equation}
\begin{equation}\label{}
\begin{aligned}
{\cal C}_L[22]:=&
\{5918918997653319693967, \\
&6686312646216567629137, \\
&7686312646216567629137, \\
&9918918997653319693967\}.
\end{aligned}
\end{equation}
\begin{equation}\label{}
\begin{aligned}
{\cal C}_L[21]:=&
\{315396334245663786197, \\
&367986315421273233617, \\
&666276812967623946997, \\
&686312646216567629137, \\
&918918997653319693967\}.
\end{aligned}
\end{equation}
It is only at the generation $20^{\text th}$ that a new contribution appears from the prime 3,
\begin{equation}\label{left-20}
\begin{aligned}
{\cal C}_L[20]:=&
\{15396334245663786197, \\
&18918997653319693967, \\
&36484957213536676883, \\
&66276812967623946997, \\
&67986315421273233617, \\
&86312646216567629137\}.
\end{aligned}
\end{equation}

\begin{figure}[t]
  \includegraphics[width=0.5\textwidth]{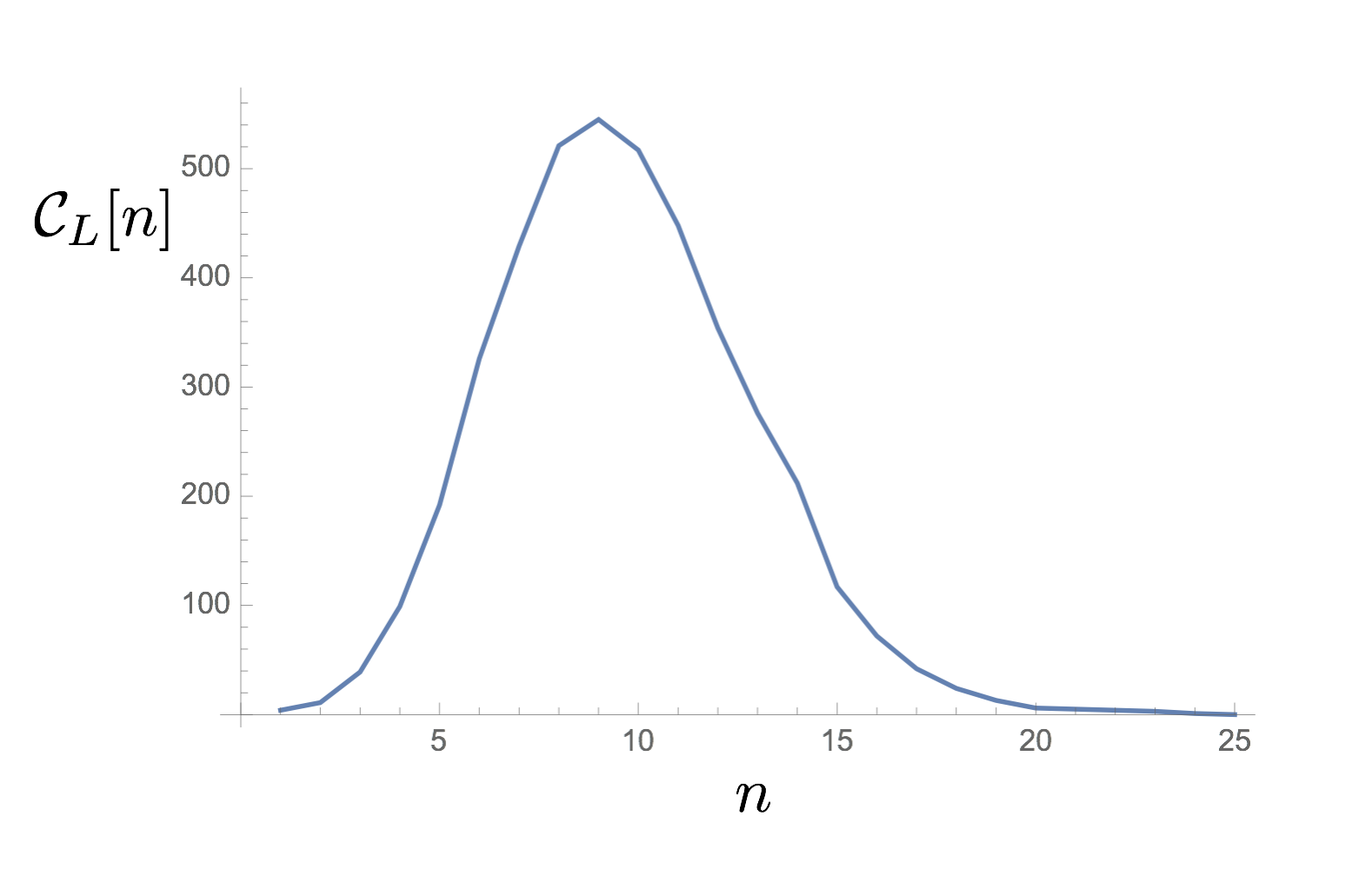}
 \caption{ The number of primes obtained by left-concatenation ${\cal C}_L[n]$ using the algorithm of Table \ref{tableL} as 
 a function of the generation number $n$. After the generation $n=24$ the left-concatenation procedure stops.}
  \label{fig:left}
\end{figure}

%%%%%%%%%%%%%%%%%%%%%%%%%%%%%%%%%%%%%%
%%%%%%%%%%%%%%%%%%%%%%%%%%%%%%%%%%%%%%
\subsection{Anomalous Left-Concatenations of Primes}
\label{sec:anomalous}
%%%%%%%%%%%%%%%%%%%%%%%%%%%%%%%%%%%%%%
%%%%%%%%%%%%%%%%%%%%%%%%%%%%%%%%%%%%%%
By definition, left-concatenations of primes that allow the 0 to be included in the set of blocks ${\cal B}_L $ are called anomalous. 
To see the effect of this inclusion of 0, let us consider some explicit example and how to make sense out of it.

Consider the 3 and left-concatenate with the 1 giving 13 that is prime. Then start introducing 0s in between 1 and 3 and test for the primality of the resulting number. As we will be soon dealing with very large number of digits, it is convenient to use some notation. Denote 0 with a subscript meaning the number of 0s coded. For instance, $\cdot 0_2\cdot = 00$, $\cdot 0_3\cdot = 000$ etc. The following series of numbers with 1 and 3 at the ends are prime numbers:
\begin{equation}\label{gap13}
\begin{aligned}
& 103 \\
&1\cdot 0_4\cdot 3 \\
&1\cdot 0_5\cdot 3 \\
&1\cdot 0_{10}\cdot 3\\
&1\cdot 0_{16}\cdot 3\\
&1\cdot 0_{17}\cdot 3\\
&1\cdot 0_{38}\cdot 3
\end{aligned}
\end{equation}
At this point, we have arrived at a prime of 40 digits. By removing the 1 on the left we have two options. Either to allow the 0s in the gap as valid degenerate digits or to discard them all as a zero to the left.
In the former case we would have a left-concatenated prime of 40 digits whereas in the latter case the prime would be only 2 digits long. This process can be repeated with the 7 replacing the 3 yielding the following series of prime numbers:
\begin{equation}\label{gap17}
\begin{aligned}
& 107 \\
&1\cdot 0_3\cdot 7 \\
&1\cdot 0_7\cdot 7 \\
&1\cdot 0_{8}\cdot 7\\
&1\cdot 0_{23}\cdot 7\\
&1\cdot 0_{59}\cdot 7\\
\end{aligned}
\end{equation}
Now, by including the 0 in ${\cal B}_L $ we can introduce a gap of 0s to the left of the left-concatenated prime number of 24 digits \eqref{largest} and produce the following prime number:
\begin{equation}\label{largest25}
1\cdot 0_{41}\cdot 357686312646216567629137
\end{equation}
If we count the 0s in the gap as degenerate digits the left-concatenated number would have 66 digits, but in the stronger case of counting all the 0s as nothing the resulting prime has still 25 digits surpassing the value of \eqref{largest}.
This process can go on to obtain new anomalous left-concatenated primes of larger number of digits. In Sec.\ref{sec:unsolved} we introduce the notion of left-chiral prime numbers and discuss its consequences.

We may call the prime numbers with a gap of 0s in between as in \eqref{gap13} and \eqref{gap17} gapped prime numbers.
Another instance of this is \eqref{largest25}. Following with this analogy with the energy spectra structures, we may think of the standard left-concatenated primes of Sec.\ref{sec:concatenatedL} as gapless prime numbers. That is, as a band of prime numbers without gaps, i.e., as a metal prime. Furthermore, we can extend the analogy by using the blocks of left-concatenated prime numbers and concatenate them to the left with gaps of zeros. For example, these are two doublets of left-concatenated primes in \eqref{left-second} separated by a gap of zeros:
\begin{equation}\label{block1}
13 \cdot 0_5 \cdot 9.
\end{equation}
We can produce more complicated structures like two prime bands separated by a gap of zeros,
\begin{equation}\label{block2}
15396334245663786197 \cdot 0_{38} \cdot 36484957213536676883,
\end{equation}
where the blocks of primes to the left and right are instances of left-concatenated prime numbers from the twentieth generation \eqref{left-20}.
A word of caution is in order: notice that in the composite prime numbers \eqref{block1} and \eqref{block2}, they are left-concatenated prime numbers by blocks but not digit by digit as in \eqref{largest25}.

By continuing enlarging the size of the gap of zeros we may increase the size of these prime numbers easily. It is tempting to
compare this sequence of prime numbers with the sequence of Mersenne prime numbers \cite{mersenne1} as an alternative route to construct the largest prime number known. Fast algorithms for searching Mersenne primes exists and with them it has been able to find the largest prime number known up to date \cite{mersenne2}.

%%%%%%%%%%%%%%%%%%%%%%%%%%%%%%%%%%%%%%
%%%%%%%%%%%%%%%%%%%%%%%%%%%%%%%%%%%%%%
\section{Unsolved Problems}
\label{sec:unsolved}
%%%%%%%%%%%%%%%%%%%%%%%%%%%%%%%%%%%%%%
%%%%%%%%%%%%%%%%%%%%%%%%%%%%%%%%%%%%%%
Chiral concatenation of prime numbers provide a method to construct sequences of prime numbers
that is amenable to simpler studies than other well-known prime number sequences. 
There are many ways to generate prime numbers sequentially. Some famous sequences are: 

\noindent i/ Dirichlet arithmetic progressions.

\noindent ii/ Fibonacci prime numbers.

\noindent iii/ Bell prime numbers.

\noindent iv/ Delannoy prime numbers.

A very difficult question is to determine whether prime sequences are infinite or not. The sequence of prime numbers itself was proven to be infinite by Euclid \cite{Euclid} with a proof that is the first notion of a theorem in the history of Mathematics.
Dirichlet arithmetic progressions extends Euclid's theorem in the sense that there are infinitely many prime numbers forming arithmetics progressions. The proof is due to Dirichlet \cite{Dirichlet} and it admits several generalizations. Fibonacci prime numbers are Fibonacci numbers that are also prime. In this case, it is not known whether there are an infinite number of Fibonacci primes \cite{Weisstein}. Similarly, Bell prime numbers are Bell numbers that are also prime. A Bell number counts the number of partitions of
a set of $n$ elements. It is also unknown whether the set of Bell prime numbers is infinite \cite{Weisstein2}. Delannoy numbers are 
generalizations of king's walks in the chess game. When they are prime they become Delannoy prime numbers and it also happens 
that it is an open question whether the Delannoy prime series is infinite \cite{Weisstein2}.

And the list of integer sequences of primes goes on an on and the question on its finiteness remains open \cite{Weisstein2}. 
In this context, what is interesting about the integer sequence of chiral concatenated prime numbers is that its finiteness 
can be analyzed exactly. Thus, the problem of finiteness for chiral concatenations of prime numbers is solved by the following results:

\noindent {\it Proposition R: Largest Right-Concatenated Prime is Finite}. 
The following statement is false: Denote $C_R[L]$ an element of the set 
${\cal C}_R[L]$, i.e., a right-concatenated prime number with $L$ digits.
$\forall n_0\in \N, \exists L\in\N$, $L>n_0$ such that if $C_R[L]$ is prime then $C_R[L+1]$
is also prime.

This is a consequence of the set ${\cal C}_R[9]=\varnothing$ being empty. There are not right-concatenated
prime numbers of arbitrary length. 
Similarly for the left-concatenated prime numbers as a consequence that the set ${\cal C}_L[25]=\varnothing$ is empty.

\noindent {\it Proposition L: Largest Left-Concatenated Prime is Finite}. 
The following statement is false: Denote $C_L[L]$ an element of the set 
${\cal C}_L[L]$, i.e., a left-concatenated prime number with $L$ digits.
$\forall n_0\in \N, \exists L\in\N$, $L>n_0$ such that if $C_L[L]$ is prime then $C_L[L+1]$
is also prime.

\noindent {\it Unsolved:} Provide a mathematical proof of Proposition R using arguments without a computer program.

\noindent {\it Unsolved:} Provide a mathematical proof of Proposition L using arguments without a computer program.

As a consequence of the notion of anomalous prime concatenations introduced in Sec.\ref{sec:anomalous} there is another interesting open question:

\noindent {\it Unsolved:} Determine whether the series of anomalous left-concatenated primes is infinite.

We may define a chirality of the prime numbers when either the sequence of right-concatenated or left-concatenated is infinite. We have seen that as both are finite, prime numbers are non-chiral. Nevetheless, there exists a clear difference in the growth of the left-concatenations and the right-concatenations of primes. In addition, we have indications to conjecture that prime numbers are left-chiral under anomalous left-concatenations.

\noindent {\it Unsolved:} Let $n_i, i=1\ldots B$ be the set of positive integers corresponding to the lengths of the blocks of left-concatenated primes chosen from the sets of generations ${\cal C}_L[n]$ with $2\leq n \leq 24$. $B$ refers to the number of bands in the analogy with energy spectra.  
Determine whether  there exist an integer $B\geq 2$ such that a $n_1-\ldots -n_B$-block left-concatenated prime number that is also left-concatenated digit by digit, with $n_i\geq 2$.

\noindent A simple example of this is the prime number \eqref{largest25} for $n_1=24, n_2=1$, but the question is for integers $n_i>1$.

It is interesting to speculate what would happen if one or both propositions L/R were false. In particular, whether the existence of infinite sequences of concatenated prime numbers would violate any of the known results in the theory of primes. This would be a way to prove propositions L/R by reductio ad absurdum. One of the best known properties of the prime numbers is the prime number theorem \cite{Gauss} that dictates asymptotically the distribution of prime numbers among the integers. This is quantified with the prime counting function denoted as $\pi(N)$ that counts how many primes there are below a given integer $N$. It behaves asymptotically as $\pi(N) \sim \frac{N}{\ln N}$ for large $N$. It turns out that the figures \ref{fig:left} and \ref{fig:right} are compatible with the prime number theorem. Let us focus on ${\cal C}_L[n]$  and the similar goes for ${\cal C}_R[n]$.
Firstly, the initial growth of ${\cal C}_L[n]$ is due to the speed proportional to the number of elements in ${\cal B}_L$. As the left-concatenation is a first order recursion, this produces an exponential growth. However, this can not last for large values of $n$ due to the prime number theorem and consequently, it must stop at some point. This is why a maximum appears in Fig.\ref{fig:left}. Secondly, after the maximum is reached, the decreasing behaviour is compatible with the asymptotic behaviour of the density of primes given by  $\frac{\pi(N)}{N} \sim \frac{1}{\ln N}$, which is a decreasing function. In particular, about $n=20$ the density of left-concatenated primes is below the density of primes.

%%%%%%%%%%%%%%%%%%%%%%%%%%%%%%%%%%%%%%
%%%%%%%%%%%%%%%%%%%%%%%%%%%%%%%%%%%%%%
\section{Conclusions}
\label{sec:conclusions}
%%%%%%%%%%%%%%%%%%%%%%%%%%%%%%%%%%%%%%
%%%%%%%%%%%%%%%%%%%%%%%%%%%%%%%%%%%%%%

The maximum generation of left-concatenated primes is 3 times longer than the right-concatenated primes: 24 vs. 8, respectively. This asymmetry is natural since the number of blocks to construct left-concatenations \eqref{left_block} from the same seed set \eqref{seed} is also bigger than the blocks used to make the right-concatenations \eqref{right_block}. None of the two concatenations produces an infinite sequence of integer sequence of prime numbers. In particular, the left-concatenation of prime numbers produces a unique prime number with the nice property that by removing its digits one by one from the left, the series of integer numbers remains primal:
357686312646216567629137. The whole series is shown in App.\ref{sec:sequence}. The limit of the 24th generation of left-concatenated primes can be surpassed by allowing for anomalous left-concatenations as introduced in Sec.\ref{sec:anomalous}. This way, an analogy with the energy spectra can be drawn in order to construct gapped primes, gapless primes, banded primes and so on and so forth. Rephrasing R.P. Feynman, there is plenty of room to the left.

%\'{\i}

\begin{acknowledgments}
I thank R. Campos, P.A.M. Casares and A. Rivas for a useful reading of the manuscript.
M.A.M.-D. acknowledges financial support from the Spanish MINECO, FIS 2017-91460-EXP, PGC2018-099169-B-I00 FIS-2018 and the CAM research consortium QUITEMAD+, Grant S2018-TCS-4243. The research of M.A.M.-D. has been supported in part by the U.S. Army Research Office through Grant No. W911N F-14-1-0103.

%\'{\i}

\end{acknowledgments}

\appendix

%%%%%%%%%%%%%%%%%%%%%%%%%%%%%%%%%%%%%%
%%%%%%%%%%%%%%%%%%%%%%%%%%%%%%%%%%%%%%
\section{List of Right-Concatenated Prime Numbers }
\label{sec:right}
%%%%%%%%%%%%%%%%%%%%%%%%%%%%%%%%%%%%%%
%%%%%%%%%%%%%%%%%%%%%%%%%%%%%%%%%%%%%%
The remaining right-concatenated prime numbers not included in Sec.\ref{sec:concatenatedR} are the following:

\begin{equation}\label{}
\begin{aligned}
{\cal C}_R[4]:=& 
\{2333, 2339, 2393, 2399, 2939, \\
&3119, 3137, 3733, 3739, 3793, 3797, \\
&5939, 7193, 7331, 7333, 7393\}.
\end{aligned}
\end{equation}

\begin{equation}\label{}
\begin{aligned}
{\cal C}_R[5]:= &
\{23333, 23339, 23399, 23993, 29399,\\
 &31193, 31379, 37337, 37339, 37397, \\
 &59393, 59399, 71933, 73331, 73939\}.
\end{aligned}
\end{equation}

\begin{equation}\label{}
\begin{aligned}
{\cal C}_R[6]:= &
\{233993, 239933, 293999, 373379, 373393, 593933, \\
&593993, 719333,  739391, 739393, 739397, 739399\}.
\end{aligned}
\end{equation}

\begin{equation}\label{}
\begin{aligned}
{\cal C}_R[7]:= &
\{2339933, 2399333, 2939999, 3733799, \\
&5939333, 7393913, 7393931,  7393933\}.
\end{aligned}
\end{equation}

\begin{equation}\label{}
\begin{aligned}
{\cal C}_R[8]:= &
\{23399339, 29399999, 37337999, \\
&59393339, 73939133\}.
\end{aligned}
\end{equation}

\vspace{90 pt}

%\vspace{5 cm}

%%%%%%%%%%%%%%%%%%%%%%%%%%%%%%%%%%%%%%
%%%%%%%%%%%%%%%%%%%%%%%%%%%%%%%%%%%%%%
\section{Largest Sequence of Left-Concatenated Primes}
\label{sec:sequence}
%%%%%%%%%%%%%%%%%%%%%%%%%%%%%%%%%%%%%%
%%%%%%%%%%%%%%%%%%%%%%%%%%%%%%%%%%%%%%

Here it is shown the resulting integer prime sequence obtained from the largest left-concatenated prime number \eqref{largest} by removing digits from the left, digit by digit.

\begin{equation}\label{sequence}
\begin{aligned}
7, \\
37,\\
137,\\
9137,\\
29137,\\
629137,\\
7629137,\\
67629137,\\
567629137,\\
6567629137,\\
16567629137,\\
216567629137,\\
6216567629137,\\
46216567629137,\\
646216567629137,\\
2646216567629137,\\
12646216567629137,\\
312646216567629137,\\
6312646216567629137,\\
86312646216567629137,\\
686312646216567629137,\\
7686312646216567629137,\\
57686312646216567629137,\\
357686312646216567629137.
\end{aligned}
\end{equation}

%%%%%%%%%%%%%%%%%%%%%%%%%%%%%%%%%%%%%%%%%%%%%%%%%%%%%%%%%%%%%%%%%%%%%%%%%%%%%%
%\begin{references}

\end{document}